\documentclass[11pt]{amsart}

\allowdisplaybreaks[4]

\newtheorem{theorem}{\bf Theorem}[section]

\newtheorem{lemma}{\bf Lemma}[section]
\newtheorem{proposition}{\bf Proposition}[section]

\theoremstyle{definition}
\newtheorem{definition}{\bf Definition}[section]

\theoremstyle{remark}

\numberwithin{equation}{section}

\newcommand{\Z}{\ensuremath{\mathbb{Z}}}
\newcommand{\C}{\ensuremath{\mathbb{C}}}
\newcommand{\R}{\ensuremath{\mathbb{R}}}

\newcommand{\N}{\ensuremath{\mathbb{N}}}
\newcommand{\Q}{\ensuremath{\mathbb{Q}}}

\def\p{\mathfrak p}
\def\P{\mathfrak P}
\def\D{\mathfrak D}

\begin{document}

\markboth{B. Behera and Q. Jahan}{Affine, quasi-affine and co-affine frames on local fields of positive characteristic}


%
%%%%%%%%%%%%%%%%%%%%%%%%%%%%%%%%%%%%%%%%%%%%%%%%%%%%%%%%%%%%%%%%%%%%

\title[AFFINE, QUASI-AFFINE AND CO-AFFINE FRAMES ON LOCAL FIELDS]{AFFINE, QUASI-AFFINE AND CO-AFFINE FRAMES ON LOCAL FIELDS OF POSITIVE CHARACTERISTIC}

\author{BISWARANJAN BEHERA}

\address{Stat-Math Unit\\ Indian Statistical Institute\\ 203 B. T. Road\\ Kolkata, 700108, India.\\
biswa@isical.ac.in}

\author{QAISER JAHAN}

\address{Department of Mathematics and Statistics\\ Indian Institute of Technology Kanpur\\ Kanpur, 208016, India.\\
qaiser@iitk.ac.in}


\begin{abstract}
The concept of quasi-affine frame in Euclidean spaces was introduced to obtain translation invariance of the discrete wavelet transform. We extend this concept to a local field $K$ of positive characteristic. We show that the affine system generated by a finite number of functions is an affine frame if and only the corresponding quasi-affine system is a quasi-affine frame. In such a case the exact frame bounds are equal. This result is obtained by using the properties of an operator associated with two such affine systems. We characterize the translation invariance of such an operator. A related concept is that of co-affine system. We show that there do not exist any co-affine frame in $L^2(K)$.
\end{abstract}
\maketitle

%
%%%%%%%%%%%%%%%%%%%%%%%%%%%%%%%%%%%%%%%%%%%%%%%%%%%%%%%%%%%%%%%%%%%%

\section{Introduction}

The concept of quasi-affine frames in $\R^n$ was introduced by Ron and Shen in~\cite{RS}, where they proved that quasi-affine frames are invariant by translations with respect to elements of $\Z^n$. They also proved that if $X$ is the affine system generated by a finite set $\Psi\subset L^2(\R^n)$ and associated with a dilation matrix $A$, and $\tilde X$ is the corresponding quasi-affine system, then $X$ is an affine frame if and only if $\tilde X$ is a quasi-affine frame, provided the Fourier transforms of the functions in $\Psi$ satisfy some mild decay conditions. Later, Chui, Shi, and St\"{o}ckler~\cite{CSS} gave an alternative proof of this fact, and more importantly, removed the assumption of the decay conditions. This result was used by Bownik in~\cite{Bow2} to provide a new characterization of multiwavelets on $L^2(\R^n)$. 

Another concept related to this theme is that of co-affine systems initially defined in~\cite{GLWW} for the case of $\R$ where the authors proved that the co-affine system can never be a frame for $L^2(\R)$. This result was subsequently extended to $L^2(\R^n)$ by Johnson~\cite{J}. Some of the other interesting articles dealing with these concepts are~\cite{Bow1, HLW, HLWW}. In this article, we extend these concepts to local fields of positive characteristic and prove analogous results. 

A field $K$ equipped with a topology is called a local field if both the additive and multiplicative groups of $K$ are locally compact abelian groups. The local fields are essentially of two types (excluding the connected local fields $\R$ and $\mathbb C$). The local fields of characteristic zero include the $p$-adic field $\Q_p$. Examples of local fields of positive characteristic are the Cantor dyadic group and the Vilenkin $p$-groups.

In order to define affine frames in a local field, we need analogous notions of translation and dilation. As we explain in the next section, we can use a prime element of a local field to serve as a dilation. In a local field of positive characteristic, we can also find a countable discrete subgroup of $K$ of the form $\{u(k):k\in\N\cup\{0\}\}$ which can be used as the translation set similar to the translation set $\Z^n$ in the standard Euclidean setup. For various aspects of wavelet analysis on a local field of positive characteristic, we refer to the articles~\cite{Beh, BJ1, BJ2, BJ3, JLJ}.

This article is organized as follows. In section~2, we provide a brief introduction to Fourier analysis on local fields. In section~3, we define affine and quasi-affine systems on a local field $K$ of positive characteristic and prove the main results of the article. We show that an affine system $X(\Psi)$ is an affine frame for $L^2(K)$ if and only if the corresponding quasi-affine system $\tilde X(\Psi)$ is a quasi-affine frame. Moreover, their exact lower and upper bounds are equal. This result also holds for Bessel families. We also characterize the translation invariance of a sesquilinear operator associated with a pair of affine systems. In section~4, we define the affine and quasi-affine dual of a finite subset $\Psi$ of $L^2(K)$ and show that a finite subset $\Phi$ of $L^2(K)$ is an affine dual of $\Psi$ if and only if it is a quasi-affine dual. In the last section, we show that $L^2(K)$ cannot have a co-affine frame.
%
%%%%%%%%%%%%%%%%%%%%%%%%%%%%%%%%%%%%%%%%%%%%%%%%%%%%%%%%%%%%%%%%%%%%

\section{Preliminaries on Local Fields}

Let $K$ be a field and a topological space. Then $K$ is called a \emph{locally compact field} or a \emph{local field} if both $K^+$ and $K^*$ are locally compact abelian groups, where $K^+$ and $K^*$ denote the additive and multiplicative groups of $K$ respectively.

If $K$ is any field and is endowed with the discrete topology, then $K$ is a local field. Further, if $K$ is connected, then $K$ is either $\R$ or $\mathbb C$. If $K$ is not connected, then it is totally disconnected. So by a local field, we mean a field $K$ which is locally compact, nondiscrete and totally disconnected.

We use the notation of the book by Taibleson~\cite{Taib}. Proofs of all the results stated in this section can be found in the books~\cite{Taib} and~\cite{RV}.

Let $K$ be a local field. Since $K^+$ is a locally compact abelian group, we choose a Haar measure $dx$ for $K^+$. If $\alpha\ne 0, \alpha\in K$, then $d(\alpha x)$ is also a Haar measure. Let $d(\alpha x)=|\alpha|dx$. We call $|\alpha|$ the \emph{absolute value} or \emph{valuation} of $\alpha$. We also let $|0| = 0$.

The map $x\rightarrow |x|$ has the following properties:
\begin{itemize}
\item[(a)] $|x| = 0$ if and only if $x=0$;
\item[(b)] $|xy|=|x||y|$ for all $x,y\in K$;
\item[(c)] $|x+y|\leq\max\{|x|, |y|\}$ for all $x,y\in K$.
\end{itemize}
Property (c) is called the \emph{ultrametric inequality}. It follows that 
\begin{equation}\label{e.max}
|x+y|=\max\{|x|, |y|\} ~\mbox{if}~|x|\ne |y|.
\end{equation}

The set $\D=\{x\in K : |x|\leq 1\}$ is called the \emph{ring of integers} in $K$. It is the unique maximal compact subring of $K$. Define $\P=\{x\in K:|x|<1\}$. The set $\P$ is called the \emph{prime ideal} in $K$. Since $K$ is totally disconnected, the set of values $|x|$ as $x$ varies over $K$ is a discrete set of the form $\{s^k: k\in\Z\}\cup \{0\}$ for some $s>0$. Hence, there is an element of $\P$ of maximal absolute value. Let $\p$ be a fixed element of maximum absolute value in $\P$. Such an element is called a \emph{prime element} of $K$. 

It can be proved that $\D$ is compact and open. Hence, $\P$ is compact and open. Therefore, the residue space $\D/\P$ is isomorphic to a finite field $GF(q)$, where $q=p^c$ for some prime $p$ and $c\in\N$. For a proof of this fact we refer to ~\cite{Taib}.

For a measurable subset $E$ of $K$, let $|E|=\int_K{\mathbf 1}_E(x)dx$, where ${\mathbf 1}_E$ is the characteristic function of $E$ and $dx$ is the Haar measure of $K$ normalized so that $|\D|=1$. Then, it is easy to see that $|\P|=q^{-1}$ and $|\p|=q^{-1}$ (see~\cite{Taib}). It follows that if $x\neq 0$, and $x\in K$, then $|x|=q^k$ for some $k\in\Z$.

Let $\D^*=\D\setminus\P=\{x\in K: |x|=1\}$. If $x\neq 0$, we can write $x=\p^k x'$, with $x'\in\D^*$. Let $\P^k=\p^k\D=\{x\in K: |x|\leq q^{-k}\}, k\in\Z$. These are called \emph{fractional ideals}. Each $\P^k$ is compact and open and is a subgroup of $K^+$ (see~\cite{RV}). 

If $K$ is a local field, then there is a nontrivial, unitary, continuous character $\chi$ on $K^+$. It can be proved that $K^+$ is self dual (see~\cite{Taib}). Let $\chi$ be a fixed character on $K^+$ that is trivial on $\D$ but is nontrivial on $\P^{-1}$. We can find such a character by starting with any nontrivial character and rescaling. We will define such a character for a local field of positive characteristic. For $y\in K$, we define $\chi_y(x)=\chi(yx)$, $x\in K$.

If $f\in L^1(K)$, then the Fourier transform of $f$ is the function $\hat f$ defined by
\[
\hat f(\xi)= \int_K f(x)\overline{\chi_{\xi}(x)}~dx.
\]
To define the Fourier transform of function in $L^2(K)$, we introduce the functions $\Phi_k$. For $k\in\Z$, let $\Phi_k$ be the characteristic function of $\mathfrak{P}^k$. For $f\in L^2(K)$, let $f_k=f\Phi_{-k}$ and define
\[
\hat{f}(\xi)=\lim\limits_{k\rightarrow\infty}\hat f_k(\xi)
=\lim\limits_{k\rightarrow \infty}\int_{\left|x\right|\leq q^k} f(x)\overline{\chi_{\xi}(x)}~d\xi,
\]
where the limit is taken in $L^2(K)$. It turns out that the Fourier transform is unitary on $L^2(K)$ (see Theorem 2.3 in~\cite{Taib}).

Let $\chi_u$ be any character on $K^+$. Since $\mathfrak{D}$ is a subgroup of $K^+$, the restriction $\chi_{u}|_\mathfrak{D}$ is a character on $\mathfrak{D}$. Also, as characters on $\mathfrak{D}, \chi_u = \chi_v$ if and only if $u-v\in \mathfrak{D}$. That is, $\chi_u=\chi_v$ if $u+\mathfrak{D}=v+\mathfrak{D}$ and $\chi_u\neq \chi_v$ if $(u+\mathfrak{D})\cap (v+\mathfrak{D})=\phi$. Hence, if $\{u(n)\}_{n=0}^{\infty}$ is a complete list of distinct coset representative of $\mathfrak{D}$ in $K^+$, then $\{\chi_{u(n)}\}_{n=0}^{\infty}$ is a list of distinct characters on $\mathfrak{D}$. It is proved in~\cite{Taib} that this list is complete. That is, we have the following proposition.

\begin{proposition}
Let $\{u(n)\}_{n=0}^{\infty}$ be a complete list of (distinct) coset representatives of $\mathfrak{D}$ in $K^+$. Then $\{\chi_{u(n)}\}_{n=0}^{\infty}$ is a complete list of (distinct) characters on $\mathfrak{D}$. Moreover, it is a complete orthonormal system on $\mathfrak{D}$.
\end{proposition}

Given such a list of characters $\{\chi_{u(n)}\}_{n=0}^{\infty}$, we define the Fourier coefficients of $f\in L^1(\mathfrak{D})$ as
\[
\hat{f}\bigl(u(n)\bigr)=\int_{\mathfrak{D}}f(x)\overline{\chi_{u(n)}(x)}dx.
\]
The series $\sum\limits_{n=0}^{\infty}\hat{f}\bigl(u(n)\bigr)\chi_{u(n)}(x)$ is called the Fourier series of $f$. From the standard $L^2$-theory for compact abelian groups we conclude that the Fourier series of $f$ converges to $f$ in $L^2(\mathfrak{D})$ and Parseval's identity holds:
\[
\int_{\mathfrak{D}}|f(x)|^2dx= \sum\limits_{n=0}^{\infty}|\hat{f}(u(n))|^2.
\]
Also, if $f\in L^1(\mathfrak{D})$ and $\hat f\bigl(u(n)\bigr)=0$ for all $n=0, 1, 2,\dots$, then $f=0$ a. e.

These results hold irrespective of the ordering of the characters. We now proceed to impose a natural order on the sequence $\{u(n)\}_{n=0}^{\infty}$. Note that $\mathfrak{D}/\mathfrak{P}$ is isomorphic to the finite field $GF(q)$ and $GF(q)$ is a $c$-dimensional vector space over the field $GF(p)$. We choose a set $\{1=\epsilon_0, \epsilon_1, \epsilon_2, \cdots, \epsilon_{c-1}\}\subset\mathfrak{D}^*$ such that span$\{\epsilon_j\}_{j=0}^{c-1}\cong GF(q)$. Let $\mathcal{U}=\{a_i: i=0, 1,\dots, q-1\}$ be any fixed full set of coset representatives of $\P$ in $\D$.
Let $\N_0=\N\cup \{0\}$. For $n\in \N_0$ such that $0\leq n< q$, we have
\[
n=a_0+a_1 p+\cdots+a_{c-1} p^{c-1},\quad 0\leq a_k<p, k=0,1,\cdots,c-1.
\]
Define
\begin{equation}\label{e.undef1}
u(n)=(a_0+a_1\epsilon_1+\cdots+a_{c-1}\epsilon_{c-1})\mathfrak{p}^{-1}.
\end{equation}
Note that $\{u(n):n=0, 1,\dots, q-1\}$ is a complete set of coset representatives of $\mathfrak{D}$ in ${\mathfrak{P}}^{-1}$. Now, for $n\geq 0$, write
\[
n=b_0+b_1q+b_2q^2+\cdots+b_sq^s,\quad 0\leq b_k<q, k=0,1,2,\cdots,s,
\]
and define
\begin{equation}\label{e.undef2}
u(n)=u(b_0)+u(b_1)\mathfrak{p}^{-1}+\cdots+u(b_s)\mathfrak{p}^{-s}.
\end{equation}

This defines $u(n)$ for all $n\in\N_0$. In general, it is not true that $u(m+n)=u(m)+u(n)$. But it follows that
\begin{equation}\label{eq.un}
u(rq^k+s)=u(r)\mathfrak{p}^{-k}+u(s)\quad{\rm if}~r\geq 0, k\geq 0~{\rm and}~0\leq s <q^k.
\end{equation}

In the following proposition we list some properties of $\{u(n): n\in\N_0\}$ which will be used later. For a proof, we refer to~\cite{BJ2}.

\begin{proposition}\label{p.un}
For $n\in\N_0$, let $u(n)$ be defined as in (\ref{e.undef1}) and (\ref{e.undef2}). Then
\begin{enumerate}
\item[{\rm(a)}] $u(n)=0$ if and only if $n=0$. If $k\geq 1$, then $|u(n)|=q^k$ if and only if $q^{k-1}\leq n < q^k$;
\item[{\rm(b)}] $\{u(k): k\in\N_0\}=\{-u(k): k\in\N_0\}$;
\item[{\rm(c)}] for a fixed $l\in\N_0$, we have $\{u(l)+u(k): k\in\N_0\}=\{u(k): k\in\N_0\}$.
\end{enumerate}
\end{proposition}

For brevity, we will write $\chi_n=\chi_{u(n)}$ for $n\in\N_0$. As mentioned before, $\{\chi_n: n\in\N_0\}$ is a complete set of characters on $\mathfrak{D}$.

Let $K$ be a local field of characteristic $p>0$ and $\epsilon_0, \epsilon_1, \dots, \epsilon_{c-1}$ be as above. We define a character $\chi$ on $K$ as follows~(see~\cite{Zheng}):
\begin{equation*}\label{chi}
\chi(\epsilon_{\mu}\mathfrak{p}^{-j})=
\left\{
\begin{array}{lll}
\exp(2\pi i/p), & \mu=0~\mbox{and}~j=1,\\
1, & \mu=1,\cdots,c-1~\mbox{or}~j\neq 1.
\end{array}
\right.
\end{equation*}
Note that $\chi$ is trivial on $\mathfrak{D}$ but nontrivial on $\mathfrak{P}^{-1}$.

In order to be able to define the concept of affine frames on local fields, we need analogous notions of translation and dilation. Since $\bigcup\limits_{j\in\Z}\mathfrak{p}^{-j}\mathfrak{D}=K$, we can regard $\mathfrak{p}^{-1}$ as the dilation (note that $|\mathfrak{p}^{-1}|=q$) and since the set $\Lambda=\{u(n): n\in\N_0\}$ is a complete list of distinct coset representatives of $\mathfrak{D}$ in $K$, it  can be treated as the translation set. Note that it follows from Proposition~\ref{p.un} that the transalation set $\Lambda$ is a subgroup of $K^+$ even though it is indexed by $\N_0$.
%
%%%%%%%%%%%%%%%%%%%%%%%%%%%%%%%%%%%%%%%%%%%%%%%%%%%%%%%%%%%%%%%%%%%%

\section{Affine Frames and Quasi-affine Frames}

For $j\in\Z$, and $y\in K$, we define the dilation operator $\delta_j$ and the translation operator $\tau_y$ on $L^2(K)$ as follows:
\[
\delta_j f(x)=q^{j/2}f(\mathfrak{p}^{-j}x) \quad\mbox{and}\quad \tau_y f(x)=f(x-y),\quad f\in L^2(K).
\]
Observe that these operators are unitary and satisfy the following commutation relation:
\begin{equation*}
\delta_j\tau_y=\tau_{\p^jy}\delta_j.
\end{equation*}
In particular, if $j<0$, then for $k\in\N_0$, we have
\begin{equation}\label{e.comm}
\delta_j\tau_{u(k)}=\tau_{u(q^{-j}k)}\delta_j.
\end{equation}
Let $f_{j,k}=\delta_j\tau_{u(k)} f$. Then
\[
f_{j,k}(x)=q^{j/2}f\bigl(\mathfrak{p}^{-j}x-u(k)\bigr), \quad j\in\Z, k\in\N_0.
\]
We also define 
\[
\tilde{f}_{j, k}=f_{j,k}=\delta_j\tau_{u(k)}f\quad{\rm if}~j\geq 0, k\in\N_0,
\]
and 
\[
\tilde{f}_{j, k}=q^{j/2}\tau_{u(k)}\delta_j f\quad{\rm if}~j<0, k\in\N_0.
\]

Let $\Psi=\{\psi^1,\psi^2,\dots,\psi^L\}$ be a finite family of functions in $L^2(K)$. The \emph{affine system} generated by $\Psi$ is the collection $X(\Psi)=\{\psi^l_{j, k}:1\leq l\leq L, j\in\Z, k\in\N_0\}$. The \emph{quasi-affine system} generated by $\Psi$ is $\tilde{X}(\Psi)=\{\tilde\psi^l_{j, k}:1\leq l\leq L, j\in\Z, k\in\N_0\}$. 

\begin{definition}
Let $\Psi\subset L^2(K)$ be a finite set. Then $X(\Psi)$ is called an \emph{affine Bessel family} if there exists a constant $B>0$ such that 
\begin{equation}\label{e.bessel}
\sum_{\eta\in X(\Psi)}|\langle f,\eta\rangle|^2\leq B\|f\|_2^2\quad{\rm for~all}~f\in L^2(K).
\end{equation}
If, in addition, there exists a constant $A>0, A\leq B$ such that 
\begin{equation}\label{e.frame}
A\|f\|_2^2\leq \sum_{\eta\in X(\Psi)}|\langle f,\eta\rangle|^2\leq B\|f\|_2^2\quad{\rm for~all}~f\in L^2(K),
\end{equation}
then $X(\Psi)$ is called an \emph{affine frame}. The largest $A$ and the smallest $B$ that can be used in the above inequalities are called the lower and upper frame bounds. The affine frame is called \emph{tight} if the lower and upper frame bounds are same. 

Similarly, $\tilde X(\Psi)$ is called a \emph{quasi-affine Bessel family} if there exists a constant $\tilde B>0$ such that~(\ref{e.bessel}) holds when $B$ is replaced by $\tilde B$ and $X(\Psi)$ is replaced by $\tilde{X}(\Psi)$. It is called a \emph{quasi-affine frame} if there exists a constant $\tilde A$ and $\tilde B>0$ such that~(\ref{e.frame}) holds when $A$ is replaced by $\tilde A$, $B$ is replaced by $\tilde B$ and $X(\Psi)$ is replaced by $\tilde{X}(\Psi)$.
\end{definition}

For two subsets $\Psi=\{\psi^1,\psi^2,\dots,\psi^L\}$ and $\Phi=\{\varphi^1, \varphi^2, \dots, \varphi^L\}$ of $L^2(K)$, we define a sesquilinear operator $K_{\Psi,\Phi}:L^2(K)\times L^2(K)\rightarrow\C$ by
\begin{equation}\label{e.kpsiphi}
K_{\Psi,\Phi}(f, g)=\sum_{l=1}^L\sum_{j\in\Z}\sum_{k\in\N_0}\langle f,\psi^l_{j,k}\rangle\langle \varphi^l_{j,k},g\rangle,\quad f, g\in L^2(K).
\end{equation}
Note that if $X(\Psi)$ and $X(\Phi)$ are affine Bessel families, then $K_{\Psi,\Phi}$ defines a bounded operator. Similarly, we define the operator $\tilde K_{\Psi,\Phi}$ by
\begin{equation}\label{e.kpsiphitilde}
\tilde K_{\Psi,\Phi}(f, g)=\sum_{l=1}^L\sum_{j\in\Z}\sum_{k\in\N_0}\langle f,\tilde\psi^l_{j,k}\rangle\langle \tilde\varphi^l_{j,k},g\rangle,\quad f, g\in L^2(K).
\end{equation}

It is easy to see that $K_{\Psi,\Phi}$ is dilation invariant, that is,  $K_{\Psi,\Phi}(\delta_N f, \delta_N g)=K_{\Psi,\Phi}(f, g)$ for all $N\in\Z$, and $\tilde K_{\Psi,\Phi}$ is invariant by translations with respect to $u(k)$, $k\in\N_0$. We write $K_{\Psi,\Psi}=K_{\Psi}$ and $\tilde K_{\Psi,\Psi}=\tilde K_{\Psi}$. 

For $j\geq 0$, let $D_j=\{q^j, q^j+1,\dots, 2q^j-1\}$. For any $j\in\Z$ and $f, g\in L^2(K)$, we define
\[
K_j(f, g)=\sum_{l=1}^L\sum_{k\in\N_0}\langle f,\psi^l_{j,k}\rangle\langle \varphi^l_{j,k},g\rangle
\]
and
\[
\tilde K_j(f, g)=\sum_{l=1}^L\sum_{k\in\N_0}\langle f,\tilde\psi^l_{j,k}\rangle\langle \tilde\varphi^l_{j,k},g\rangle.
\]

We first prove two crucial lemmas before we state and prove the main results of this article. 

\begin{lemma}\label{l.1}
Let $\Psi=\{\psi^1,\psi^2,\dots,\psi^L\}$ and $\Phi=\{\varphi^1, \varphi^2, \dots, \varphi^L\}$ be two subsets of $L^2(K)$. Fix $J\in\N$. Then for all $j\geq-J$ and $f, g\in L^2(K)$, we have
\[
\tilde K_j(f, g)=q^{-J}\sum_{\nu\in D_J}K_j(\tau_{u(\nu)}f, \tau_{u(\nu)}g).
\]
\end{lemma}

\begin{proof}
For $j\geq 0$, $\tilde K_j(f, g)=K_j(f, g)=K_j(\tau_{u(\nu)}f, \tau_{u(\nu)}g)$ for any $\nu\in\N_0$. Now, for any integer $j$ such that $-J\leq j<0$, $K_j$ is invariant with respect to translation by $u(q^{-j}\nu)=\p^ju(\nu)$, $\nu\in\N_0$. That is,
\[
K_j(\tau_{\p^ju(\nu)}f, \tau_{\p^ju(\nu)}g)=K_j(f, g),\quad \nu\in\N_0.
\]

Note that, for any $m\leq J$, 
\[
\sum_{\nu=0}^{q^J-1}a_\nu=\sum_{\lambda=0}^{q^m-1}\sum_{\mu=0}^{q^{J-m}-1}a_{\mu q^m+\lambda}.
\]
Hence,
\[
\sum_{\nu\in D_J}a_\nu=\sum_{\nu=q^J}^{2q^J-1}a_\nu=\sum_{\lambda=0}^{q^m-1}\sum_{\mu=0}^{q^{J-m}-1}a_{\mu q^m+\lambda+q^J}.
\]
Therefore, we have
\begin{eqnarray*}
q^{-J}\sum_{\nu\in D_J}K_j(\tau_{u(\nu)}f, \tau_{u(\nu)}g) 
& = & q^{-J}\sum_{\lambda=0}^{q^{-j}-1}\sum_{\mu=0}^{q^{J+j}-1}K_j(\tau_{u(\mu q^{-j}+\lambda+q^J)}f, \tau_{u(\mu q^{-j}+\lambda+q^J)}g) \\
& = & q^{-J}\sum_{\lambda=q^{-j}}^{2q^{-j}-1}\sum_{\mu=q^{J+j}-1}^{2q^{J+j}-2}K_j(\tau_{u(\mu q^{-j}+\lambda)}f, \tau_{u(\mu q^{-j}+\lambda)}g).
\end{eqnarray*}

Since $\lambda\in\{q^{-j}, q^{-j}+1, \dots, 2q^{-j}-1\}$ and $\mu\in\{q^{J+j}-1, q^{J+j}, 2q^{J+j}-2\}$, we have $\mu=q^{J+j}+r$ and $\lambda=q^{-j}+s$, where $r\in\{-1, 0, 1, \dots, q^{J+j}-2\}$ and $s\in\{0, 1,\dots, q^{-j}-1\}$. Hence, $\mu q^{-j}+\lambda=q^J+rq^{-j}+q^{-j}+s=(q^{J+j}+r+1)q^{-j}+s$ so that $u(\mu q^{-j}+\lambda)=u(q^{J+j}+r+1)\p^j+u(\lambda-q^{-j})$, by~(\ref{eq.un}). Therefore, 
\begin{eqnarray*}
q^{-J}\sum_{\nu\in D_J}K_j(\tau_{u(\nu)}f, \tau_{u(\nu)}g) & = & q^{-J}\sum_{\lambda=q^{-j}}^{2q^{-j}-1}\sum_{\mu=q^{J+j}-1}^{2q^{J+j}-2}K_j(\tau_{u(\lambda-q^{-j})}f, \tau_{u(\lambda-q^{-j})}g) \\
& = & q^{j}\sum_{\lambda=0}^{q^{-j}-1}K_j(\tau_{u(\lambda)}f, \tau_{u(\lambda)}g) \\
& = & q^{j}\sum_{\lambda=0}^{q^{-j}-1}K_j(\tau_{-u(\lambda)}f, \tau_{-u(\lambda)}g) \\
& = & q^{j}\sum_{\lambda=0}^{q^{-j}-1}\sum_{l=1}^L\sum_{k\in\N_0}\langle\tau_{-u(\lambda)} f,\psi^l_{j,k}\rangle\langle \varphi^l_{j,k},\tau_{-u(\lambda)}g\rangle \\
& = & \sum_{\lambda=0}^{q^{-j}-1}\sum_{l=1}^L\sum_{k\in\N_0}\langle f,\tilde\psi^l_{j,q^{-j}k+\lambda}\rangle\langle \tilde\varphi^l_{j,q^{-j}k+\lambda}, g\rangle \\
& = & \sum_{l=1}^L\sum_{k\in\N_0}\langle f,\tilde\psi^l_{j,k}\rangle\langle \tilde\varphi^l_{j,k}, g\rangle = \tilde K_j(f, g).
\end{eqnarray*}
We have used the following two facts in the series of equalities above:
\begin{enumerate}
\item[(i)]  $\{u(\lambda):0\leq\lambda\leq q^m-1\}=\{-u(\lambda):0\leq\lambda\leq q^m-1\}$ for any $m\in\N$.
\item[(ii)] $q^{j/2}\langle\tau_{-u(\lambda)} f,\psi^l_{j,k}\rangle=\langle f,\tilde\psi^l_{j,q^{-j}k+\lambda}\rangle$.
\end{enumerate}
This completes the proof of the lemma. 
\end{proof} 

In the following lemma we prove two important properties of the operators $K_j$ and $\tilde K_j$ when $\Phi=\Psi$.

\begin{lemma}\label{l.2}
Let $\Psi=\{\psi^1,\psi^2,\dots,\psi^L\}\subset L^2(K)$. Put $\Phi=\Psi$ in the definitions of $K_j$ and $\tilde K_j$. If $f\in L^2(K)$ has compact support, then
\begin{enumerate}
\item[(a)] $\lim\limits_{N\rightarrow\infty}\sum\limits_{j<0}\tilde K_j(\delta_N f, \delta_N f)=0$.
\item[(b)] $\lim\limits_{N\rightarrow\infty} q^{-N}\sum\limits_{j<-N}\sum\limits_{\nu\in D_N}K_j(\tau_{u(\nu)} f, \tau_{u(\nu)} f)=0$.
\end{enumerate}
\end{lemma}

\begin{proof}
We have,
\[
\sum\limits_{j<0}\tilde K_j(\delta_N f, \delta_N f)=\sum\limits_{j<0}\sum_{l=1}^L\sum_{k\in\N_0}|\langle \delta_N f, q^{j/2}\tau_{u(k)}\delta_j\psi^l\rangle|^2.
\]
Let $\Omega={\rm supp}~f$. Since
\[
\langle \delta_N f, \tau_{u(k)}\delta_j\psi^l\rangle 
= \langle f, \delta_{-N}\tau_{u(k)}\delta_j\psi^l\rangle 
= \langle f, \tau_{\p^{-N}u(k)}\delta_{-N}\delta_j\psi^l\rangle 
= \langle\tau_{-\p^{-N}u(k)}f, \delta_{j-N}\psi^l\rangle,
\]
we have,
\begin{eqnarray*}
|\langle \delta_N f, \tau_{u(k)}\delta_j\psi^l\rangle|^2
&   =  & \Bigl|\int_K f(x+\p^{-N}u(k))\overline{\delta_{j-N}\psi^l(x)}dx\Bigr|^2 \\
& \leq & \|f\|_2^2\int_{\Omega-\p^{-N}u(k)}|\delta_{j-N}\psi^l(x)|^2 dx \\
&   =  & \|f\|_2^2\int_{\Omega-\p^{-N}u(k)}|q^{(j-N)/2}\psi^l(\p^{-(j-N)x})|^2 dx \\
&   =  & \|f\|_2^2\int_{\p^{-(j-N)}(\Omega-\p^{-N}u(k))}|\psi^l(x)|^2 dx \\
&   =  & \|f\|_2^2\int_{\p^{-(j-N)}\Omega-\p^{-j}u(k)}|\psi^l(x)|^2 dx. 
\end{eqnarray*}
Thus,
\[
\sum\limits_{j<0}\tilde K_j(\delta_N f, \delta_N f) \leq \|f\|_2^2\sum\limits_{j<0}q^j\sum_{k\in\N_0}\int_{\p^{-(j-N)}\Omega-\p^{-j}u(k)}\sum_{l=1}^L|\psi^l(x)|^2 dx.
\]

Note that $\p^{-(j-N)}\Omega-\p^{-j}u(k)=\p^{-j}(\p^N\Omega-u(k))$. Since $\Omega$ is compact and $|\p^N\Omega|=q^{-N}|\Omega|$, we can choose $N_0$ large enough so that $\p^N\Omega\subseteq\D$ if $N\geq N_0$. Since $\{\D+u(k):k\in\N_0\}$ is a disjoint collection, it follows that $\{\p^{-j}(\p^N\Omega-u(k)):k\in\N_0\}$ is also a disjoint collection. Hence,
\begin{eqnarray*}
\sum\limits_{j<0}\tilde K_j(\delta_N f, \delta_N f)
& \leq & \|f\|_2^2\sum\limits_{j<0}q^j\int_{\bigcup\limits_{k\in\N_0}\p^{-j}(\p^N\Omega-u(k))}\sum_{l=1}^L|\psi^l(x)|^2 dx \\
&  =   & \|f\|_2^2\int_KF_N(x)\sum_{l=1}^L|\psi^l(x)|^2 dx,
\end{eqnarray*}
where
\[
F_N=\sum\limits_{j<0}q^j{\mathbf 1}_{\bigcup\limits_{k\in\N_0}\p^{-j}(\p^N\Omega-u(k))}.
\]

Observe that $|F_N(x)|\leq \sum\limits_{j<0}q^j=\frac{1}{q-1}$. Since $\sum\limits_{l=1}^L|\psi^l|^2\in L^1(K)$, if we can show that $F_N\rightarrow 0$ a.e. as $N\rightarrow\infty$, then by Lebesgue dominated convergence theorem, the last integral above will converge to $0$ as $N\rightarrow\infty$.

Let $E=\{x\in K: x=-\p^{-j}u(k)~{\rm for~some}~j<0~{\rm and}~k\in\N_0\}$. If $x\not\in E$, then $\p^j x+u(k)\not=0$ for any $j<0$ and $k\in\N_0$ so that $|\p^j x+u(k)|=q^r$ for some $r\in\Z$. Thus, $\p^j x+u(k)\not\in\p^N\Omega$ if $N>-r$. That is,  $x\not\in \p^{-j}(\p^N\Omega-u(k))$ if $N>-r$. Since $E$ is a set of measure zero, it follows that $F_N\rightarrow 0$ a.e. as $N\rightarrow\infty$. This proves part (a) of the lemma.

To prove part (b), we observe that 
\begin{eqnarray*}
q^{-N}\sum\limits_{j<-N}\sum\limits_{\nu\in D_N}K_j(\tau_{u(\nu)} f, \tau_{u(\nu)} f) 
& = & q^{-N}\sum\limits_{\nu=q^N}^{2q^N-1}\sum\limits_{j<-N}\sum_{l=1}^L\sum_{k\in\N_0}|\langle\tau_{u(\nu)}f,\psi^l_{j,k}\rangle|^2.
\end{eqnarray*}
An easy calculation as in part (a) gives us
\[
\langle\tau_{u(\nu)}f,\psi^l_{j,k}\rangle = \langle\tau_{u(\nu)-p^ju(k)}f, \delta_j\psi^l\rangle,
\]
so that
\begin{eqnarray*}
|\langle\tau_{u(\nu)}f,\psi^l_{j,k}\rangle|^2
&   =   & |\langle\tau_{u(\nu)-p^ju(k)}f, \delta_j\psi^l\rangle|^2 \\
& \leq  & \|f\|_2^2\int_{\Omega+u(\nu)-\p^ju(k)}|\delta_j\psi^l(x)|^2 dx \\
&  =    & \|f\|_2^2\int_{\p^{-j}(\Omega+u(\nu)-\p^ju(k))}|\psi^l(x)|^2 dx.
\end{eqnarray*}
Hence,
\begin{eqnarray*}
\lefteqn{
q^{-N}\sum\limits_{j<-N}\sum\limits_{\nu\in D_N}K_j(\tau_{u(\nu)} f, \tau_{u(\nu)} f)} \\
& \leq & q^{-N}\|f\|_2^2\sum\limits_{\nu=q^N}^{2q^N-1}\sum\limits_{j<-N}\sum_{k\in\N_0}\int_{\p^{-j}(\Omega+u(\nu)-\p^ju(k))}\sum_{l=1}^L|\psi^l(x)|^2 dx\\
&   =  & \|f\|_2^2\int_K G_N(x)\sum_{l=1}^L|\psi^l(x)|^2 dx,
\end{eqnarray*}
where
\[
G_N=q^{-N}\sum\limits_{\nu=q^N}^{2q^N-1}\sum\limits_{j<-N}\sum_{k\in\N_0}{\mathbf 1}_{\p^{-j}(\Omega+u(\nu))-u(k)}.
\]
As in part (a), to complete the proof, we need to show that $G_N\rightarrow 0$ a.e. as $N\rightarrow\infty$. 

Note that $\p^N\Omega\subseteq\D$ if $N\geq N_0$. For such an $N$, consider the set $\p^N(\Omega+u(\nu))$, where $\nu\in D_N$. If $x\in\p^N(\Omega+u(\nu))$, then $x=y+\p^Nu(\nu)$ for some $y\in\p^N\Omega\subseteq\D$. Since $|y|\leq 1$ and $|\p^Nu(\nu)|=q^{-N}\cdot q^{N+1}=q$, we have $|x|=q$, by~(\ref{e.max}). Thus, $\p^N(\Omega+u(\nu))\subset\P^{-1}\setminus\D=\p^{-1}\D\setminus\D=\p^{-1}\D^*$ so that $\p^{-j}(\Omega+u(\nu))\subseteq\p^{-j-N-1}\D^*$ for any $j<N$.

For $N\geq N_0$, fix $j<-N$ and $k\in\N_0$. Then
\[
\p^{-j}(\Omega+u(\nu))-u(k)\subseteq \p^{-j}(\p^{-N_0}\D+u(\nu))-u(k).
\]
Note that since $\Omega$ is compact, each $\Omega+u(k_0)$ can intersect with only finitely many sets of the form $\Omega+u(k)$, $k\in\N_0$. So there exists an integer $d\in\N$ such that each $x\in K$ can belong to at most $d$ such sets. Thus, in particular, any $x\in K$ can belong to at most $d$ sets in the collection $\{\p^{-j}(\Omega+u(\nu))-u(k):\nu\in D_N\}$. Each of these sets is contained in $\p^{-j-N-1}\D^*-u(k)$, hence so is their union. Thus,
\[
\sum\limits_{\nu=q^N}^{2q^N-1}{\mathbf 1}_{\p^{-j}(\Omega+u(\nu))-u(k)}\leq d{\mathbf 1}_{\p^{-j-N-1}\D^*-u(k)}.
\]
Now, 
\begin{eqnarray*}
\sum\limits_{j<-N}{\mathbf 1}_{\p^{-j-N-1}\D^*-u(k)} 
& = & \sum\limits_{i>0}{\mathbf 1}_{\p^{i-1}\D^*-u(k)} \\
& = & {\mathbf 1}_{\bigcup\limits_{i>0}\p^{i-1}\D^*-u(k)} = {\mathbf 1}_{\D^*-u(k)}\leq 1,
\end{eqnarray*}
since $\{\D+u(k):k\in\N_0\}$ is a partition of $K$ and $\D^*\subset\D$.

Collecting all these estimates, we get
\[
G_N(x)\leq q^{-N}\cdot d.
\]
Therefore, $G_N(x)\rightarrow 0$ as $N\rightarrow\infty$ uniformly in $x$. This completes the proof of the lemma.
\end{proof}

The following theorem shows the relationship between affine and quasi-affine frames.
\begin{theorem}\label{t.aframe}
Let $\Psi=\{\psi^1,\psi^2,\dots,\psi^L\}\subset L^2(K)$. Then
\begin{enumerate}
\item[(a)] $X(\Psi)$ is a Bessel family if and only if $\tilde X(\Psi)$ is a Bessel family. Moreover, their exact upper bounds are equal.
\item[(b)] $X(\Psi)$ is an affine frame if and only if $\tilde X(\Psi)$ is a quasi-affine frame. Moreover, their exact lower and upper bounds are equal.
\end{enumerate}
\end{theorem}

\begin{proof}
Put $\Phi=\Psi$ in the definitions of $K_j$ and $\tilde K_j$. Suppose that $X(\Psi)$ is a Bessel family with upper bound $B\geq 0$. Then, by Lemma~\ref{l.1}, for all $f\in L^2(K)$, we have
\begin{eqnarray*}
\tilde K_\Psi(f, f)
& = & \sum_{j=-\infty}^\infty\tilde K_j(f, f)=\lim_{J\rightarrow\infty}\sum_{j\geq -J}\tilde K_j(f, f)\\
& = & \lim_{J\rightarrow\infty}q^{-J}\sum_{\nu\in D_J}\sum_{j\geq -J} K_j(\tau_{u(\nu)}f, \tau_{u(\nu)}f)\\
& \leq & \lim_{J\rightarrow\infty}q^{-J}\sum_{\nu\in D_J}K_\Psi(\tau_{u(\nu)}f, \tau_{u(\nu)}f)\\
& \leq & \lim_{J\rightarrow\infty}q^{-J}\sum_{\nu\in D_J}B\|\tau_{u(\nu)}f\|_2^2=B\|f\|_2^2.
\end{eqnarray*}
Thus, the quasi-affine system $\tilde X(\Psi)$ is also a Bessel family with upper bound $B$.

Conversely, let us assume that $\tilde X(\Psi)$ is a Bessel family with upper bound $C\geq 0$. Further, assume that there exists $f\in L^2(K)$ with $\|f\|_2=1$ and $K_\Psi(f, f)>C$. We will get a contradiction. We have
\begin{eqnarray*}
\sum_{j=-N}^\infty K_j(f, f)
& = & \sum_{j=-N}^\infty\sum_{l=1}^L\sum_{k\in\N_0}|\langle f,\psi^l_{j,k}\rangle|^2 \\
& = & \sum_{j=0}^\infty\sum_{l=1}^L\sum_{k\in\N_0}|\langle f,\psi^l_{j-N,k}\rangle|^2 \\
& = & \sum_{j=0}^\infty\sum_{l=1}^L\sum_{k\in\N_0}|\langle\delta_N f, \psi^l_{j,k}\rangle|^2 \\
& = & \sum_{j=0}^\infty K_j(\delta_N f, \delta_N f).
\end{eqnarray*}
Since $K_\Psi(f, f)=\lim\limits_{N\rightarrow\infty}\sum\limits_{j=-N}^\infty K_j(f, f)>C$, there exists $N\in\N$ such that 
\[
\sum\limits_{j=-N}^\infty K_j(f, f)=\sum\limits_{j=0}^\infty K_j(\delta_N f, \delta_N f)>C.
\]
Now,
\[
\tilde K_\Psi(\delta_N f, \delta_N f)
\geq \sum_{j=0}^\infty\tilde K_j(\delta_N f, \delta_N f) 
=  \sum_{j=0}^\infty K_j(\delta_N f, \delta_N f) >C.
\]

If $g=\delta_Nf$, then we have $\|g\|_2=\|\delta_Nf\|_2=\|f\|_2=1$ but $\tilde K_\Psi(g, g)>C$. This is a contradiction to the fact that $\tilde X(\Psi)$ is a Bessel family with upper bound $C$. This proves part (a) of the theorem.

We will now prove part (b). We have dealt with the upper bounds in part (a). So we need only to consider the lower bounds $\tilde A$ and $A$. Suppose that $X(\Psi)$ is an affine frame with lower frame bound $A$. Then, for all $f\in L^2(K)$ with compact support, we have
\begin{eqnarray*}
\tilde K_\Psi(f, f)
& = & \lim_{J\rightarrow\infty}q^{-J}\sum_{\nu\in D_J}\sum_{j\geq -J} K_j(\tau_{u(\nu)}f, \tau_{u(\nu)}f)\\
& = & \lim_{J\rightarrow\infty}q^{-J}\sum_{\nu\in D_J}\sum_{j\in\Z} K_j(\tau_{u(\nu)}f, \tau_{u(\nu)}f) \quad{\rm(by~Lemma~\ref{l.2}(b))}\\
& = & \lim_{J\rightarrow\infty}q^{-J}\sum_{\nu\in D_J} K_\Psi(\tau_{u(\nu)}f, \tau_{u(\nu)}f)\\
& \geq & \lim_{J\rightarrow\infty}q^{-J}\sum_{\nu\in D_J}A\|\tau_{u(\nu)}f\|_2^2=A\|f\|_2^2.
\end{eqnarray*}
The set of all such $f$ is dense in $L^2(K)$. So this holds for all $f\in L^2(K)$. Hence, $\tilde A\geq A$. 

To show that $\tilde A\leq A$, we assume that it is not true and get a contradiction. Thus, there exists $\epsilon>0$, $f\in L^2(K)$ with $\|f\|=1$ such that  
\[
K_\Psi(f, f)\leq \tilde A-\epsilon.
\]
Without loss of generality, we can assume that $f$ has compact support (otherwise, for any compact set $\Omega$, we consider $f{\mathbf 1}_\Omega)$. Since $K_\Psi$ is dilation invariant, we also get
\[
K_\Psi(\delta_Nf, \delta_Nf)\leq \tilde A-\epsilon.
\]
By Lemma~\ref{l.2}(a), there exists $N\in\N$ such that $\sum\limits_{j<0}\tilde K_j(\delta_Nf, \delta_Nf)<\epsilon/2$. Hence,
\begin{eqnarray*}
\tilde K_\Psi(\delta_Nf, \delta_Nf)
& < & \sum_{j\geq 0}\tilde Kj(\delta_Nf, \delta_Nf)+\epsilon/2 \\
& = & \sum_{j\geq 0}Kj(\delta_Nf, \delta_Nf)+\epsilon/2 \\
& \leq & K_\Psi(\delta_Nf, \delta_Nf)+\epsilon/2 \leq \tilde A-\epsilon/2.
\end{eqnarray*}
This contradicts the definition of the lower bound $\tilde A$ of $\tilde X(\Psi)$ and completes the proof of the theorm.
\end{proof}

We have observed earlier that $K_{\Psi,\Phi}$ is dilation invariant and $\tilde K_{\Psi,\Phi}$ is invariant by translations with respect to $u(k)$, $k\in\N_0$. In the next theorem, we show that a necessary and sufficient condition for the translation invariance of $K_{\Psi,\Phi}$ is that the operators $K_{\Psi,\Phi}$ and $\tilde K_{\Psi,\Phi}$ coincide.

\begin{theorem}\label{t.transinv}
Let $\Psi=\{\psi^1,\psi^2,\dots,\psi^L\}$ and $\Phi=\{\varphi^1,\varphi^2,\dots,\varphi^L\}$ generate two affine Bessel families. Then $K_{\Psi,\Phi}$ is translation invariant if and only if $K_{\Psi,\Phi}=\tilde K_{\Psi,\Phi}$.
\end{theorem}
\begin{proof}
Suppose that $K_{\Psi,\Phi}$ is translation invariant. Then, as in the proof of Theorem~\ref{t.aframe}, for all $f, g\in L^2(K)$ with compact support, we have
\begin{eqnarray*}
\tilde K_{\Psi,\Phi}(f, g)
& = & \lim_{J\rightarrow\infty}q^{-J}\sum_{\nu\in D_J} K_{\Psi,\Phi}(\tau_{u(\nu)}f, \tau_{u(\nu)}g)\\
& = & \lim_{J\rightarrow\infty}q^{-J}\sum_{\nu\in D_J} K_{\Psi,\Phi}(f, g)=K_{\Psi,\Phi}(f, g),
\end{eqnarray*}
where we have used the translation invariance of $K_{\Psi,\Phi}$. By density and the boundedness of the operators $K_{\Psi,\Phi}$ and $\tilde K_{\Psi,\Phi}$, the equality holds for all $f, g\in L^2(K)$.

Conversely, assume that $K_{\Psi,\Phi}=\tilde K_{\Psi,\Phi}$. Then for $m\in\N_0$,
\begin{eqnarray*}
K_{\Psi,\Phi}(\tau_{u(m)}f, \tau_{u(m)}g)
& = & \tilde K_{\Psi,\Phi}(\tau_{u(m)}f, \tau_{u(m)}g) \\
& = & \sum_{l=1}^L\sum_{j\geq 0}\sum_{k\in\N_0}\langle \tau_{u(m)}f, \delta_j\tau_{u(k)}\psi^l\rangle \langle\delta_j\tau_{u(k)}\varphi^l, \tau_{u(m)}g\rangle \\
&   & + \sum_{l=1}^L\sum_{j<0}\sum_{k\in\N_0}\langle \tau_{u(m)}f, q^{j/2}\tau_{u(k)}\delta_j\psi^l\rangle \langle q^{j/2}\tau_{u(k)}\delta_j\varphi^l, \tau_{u(m)}g\rangle. 
\end{eqnarray*}
Since $m\in N_0$, by Proposition~\ref{p.un}(b), there exists a unique $m_0\in\N_0$ such that $-u(m)=u(m_0)$. Hence, in the first sum, we have 
\begin{eqnarray*}
\langle \tau_{u(m)}f, \delta_j\tau_{u(k)}\psi^l\rangle
& = & \langle f, \tau_{u(m_0)}\delta_j\tau_{u(k)}\psi^l\rangle \\
& = & \langle f, \delta_j\tau_{u(q^jm_0)+u(k)}\psi^l\rangle.
\end{eqnarray*}
Similarly, $\langle\delta_j\tau_{u(k)}\varphi^l, \tau_{u(m)}g\rangle=\langle\delta_j\tau_{u(q^jm_0)+u(k)}\varphi^l, g\rangle$.

For a fixed $j\geq 0$, we have $\{u(k)+u(q^jm_0):k\in\N_0\}= \{u(k):k\in\N_0\}$, by Proposition~\ref{p.un}(c). Hence, for a fixed $j\geq 0$, we have
\begin{eqnarray*}
\sum_{k\in\N_0}\langle \tau_{u(m)}f, \delta_j\tau_{u(k)}\psi^l\rangle \langle\delta_j\tau_{u(k)}\varphi^l, \tau_{u(m)}g\rangle
& = & \sum_{k\in\N_0}\langle f, \delta_j\tau_{u(k)}\psi^l\rangle \langle\delta_j\tau_{u(k)}\varphi^l, g\rangle \\
& = & \sum_{k\in\N_0}\langle f, \psi^l_{j, k}\rangle \langle\varphi^l_{j, k}, g\rangle.
\end{eqnarray*}

In the second sum, we have
\[
\langle \tau_{u(m)}f, q^{j/2}\tau_{u(k)}\delta_j\psi^l\rangle=\langle f, q^{j/2}\tau_{u(m_0)+u(k)}\delta_j\psi^l\rangle.
\]
By a similar argument as above, we get for each $j<0$,
\[
\sum_{k\in\N_0}\langle \tau_{u(m)}f, q^{j/2}\tau_{u(k)}\delta_j\psi^l\rangle \langle q^{j/2}\tau_{u(k)}\delta_j\varphi^l, \tau_{u(m)}g\rangle
=\sum_{k\in\N_0}\langle f, \tilde\psi^l_{j, k}\rangle \langle\tilde\varphi^l_{j, k}, g\rangle.
\]
Hence, we have 
\[
K_{\Psi,\Phi}(\tau_{u(m)}f, \tau_{u(m)}g)=\tilde K_{\Psi,\Phi}(f, g)=K_{\Psi,\Phi}(f, g).
\]

This proves that $K_{\Psi,\Phi}$ is invariant by translations with respect to $u(m)$, where $m\in\N_0$. Since $K_{\Psi,\Phi}$ is invariant with respect to dilations, it follows that it is invariant with respect to all $x$ of the form $x=\p^ju(m)$, $m\in\N_0$, $j\in\Z$. But such elements are dense in $K$. This can be seen as follows. Since $\{\D+u(m):m\in\N_0\}$ is a partition of $K$, $\{\p^j\D+\p^ju(m):m\in\N_0\}$ is also a partition of $K$ for any $j\in\Z$. Hence, if $x\in K$, then for each $j\in\Z$, there exists a unique $m\in\N_0$ and $y\in\D$ such that $x=\p^jy+\p^ju(m)$ so that $|x-\p^ju(m)|=|\p^jy|=q^{-j}|y|$. Since $|y|\leq 1$, we can choose $j$ sufficiently large to make $|x-\p^ju(m)|$ as small as we want.

Now, since $K_{\Psi,\Phi}$ is a bounded operator and translation is a continuous operation, it follows that $K_{\Psi,\Phi}$ is invariant with respect to translation by all elements of $K$. This completes the proof of the theorem.
\end{proof}
%
%%%%%%%%%%%%%%%%%%%%%%%%%%%%%%%%%%%%%%%%%%%%%%%%%%%%%%%%%%%%%%%%%%%%

\section{Affine and Quasi-affine Duals}

In this section, we define the affine dual and quasi-affine dual of a finite subset $\Psi$ of $L^2(K)$ and show that a finite subset $\Phi$ of $L^2(K)$ with cardinality same as that of $\Psi$ is an affine dual of $\Psi$ if and only if it is a quasi-affine dual of $\Psi$. 
\begin{definition}
Let $\Psi=\{\psi^1,\psi^2,\dots,\psi^L\}$ and $\Phi=\{\varphi^1, \varphi^2, \dots, \varphi^L\}$ be two subsets of $L^2(K)$ such that $X(\Psi)$ and $X(\Phi)$ are Bessel families. Then $\Phi$ is called an \emph{affine dual} of $\Psi$ if $K_{\Psi,\Phi}(f, g)=\langle f, g\rangle$ for all $f, g\in L^2(K)$, that is,
\begin{equation}\label{e.adual}
\sum_{l=1}^L\sum_{j\in\Z}\sum_{k\in\N_0}\langle f, \psi^l_{j,k}\rangle \langle\varphi^l_{j,k}, g\rangle\quad{\rm for~all~}f, g\in L^2(K).
\end{equation}
We say that $\Phi$ is a \emph{quasi-affine dual} of $\Psi$ if $\tilde K_{\Psi,\Phi}(f, g)=\langle f, g\rangle$ for all $f, g\in L^2(K)$, that is, 
\begin{equation}\label{e.qadual}
\sum_{l=1}^L\sum_{j\in\Z}\sum_{k\in\N_0}\langle f, \tilde\psi^l_{j,k}\rangle \langle\tilde\varphi^l_{j,k}, g\rangle\quad{\rm for~all~}f, g\in L^2(K).
\end{equation}
\end{definition}

Since $K_{\Psi,\Phi}$ and $\tilde K_{\Psi,\Phi}$ are sesquilinear operators, it follows from the polarization identity that~(\ref{e.adual}) or~(\ref{e.qadual}) holds if and only if it holds for all $f=g$ in $L^2(K)$.

\begin{theorem}
Let $\Psi=\{\psi^1,\psi^2,\dots,\psi^L\}\subset L^2(K)$ generate an affine Bessel family. Then $\Phi=\{\varphi^1, \varphi^2, \dots, \varphi^L\}$ is an affine dual of $\Psi$ if and only if it is a quasi-affine dual of $\Psi$.
\end{theorem}

\begin{proof} 
We first assume that $\Phi$ is an affine dual of $\Psi$. So $K_{\Psi,\Phi}(f, g)=\langle f, g\rangle$ for all $f, g\in L^2(K)$. Since $\langle\tau_y f, \tau_y g\rangle=\langle f, g\rangle$ for all $y\in K$ and for all $f, g\in L^2(K)$, it follows that $K_{\Psi,\Phi}$ is translation invariant. Hence, by Theorem~\ref{t.transinv}, we have
\[
\tilde K_{\Psi,\Phi}(f, g)=K_{\Psi,\Phi}(f, g)=\langle f, g\rangle \quad{\rm for~all~}f, g\in L^2(K).
\]
Therefore, $\Phi$ is a quasi-affine dual of $\Psi$.

Conversely, assume that $\Phi$ is a quasi-affine dual of $\Psi$. Let $f\in L^2(K)$ be a function with compact support. By Lemma~\ref{l.2}(a), we have 
\[
\sum_{j<0}\tilde Kj(\delta_N f, \delta_N f)\rightarrow 0~{\rm as}~N\rightarrow\infty.
\]
That is,
\begin{equation}\label{e.dual}
\sum_{j<0}\sum_{l=1}^L\sum_{k\in\N_0}\langle \delta_N f, \tilde\psi^l_{j,k}\rangle \langle\tilde\varphi^l_{j,k}, \delta_N f\rangle\rightarrow 0~{\rm as}~N\rightarrow\infty.
\end{equation}
Now, since $\Phi$ is a quasi-affine dual of $\Psi$, we have
\begin{eqnarray*}
&  & \|f\|^2=\|\delta_N f\|^2=\langle\delta_N f, \delta_N f\rangle=\sum_{l=1}^L\sum_{j\in\Z}\sum_{k\in\N_0}\langle \delta_N f, \tilde\psi^l_{j,k}\rangle \langle\tilde\varphi^l_{j,k}, \delta_N f\rangle \\
& = & \sum_{l=1}^L\sum_{j\geq 0}\sum_{k\in\N_0}\langle \delta_N f, \tilde\psi^l_{j,k}\rangle \langle\tilde\varphi^l_{j,k}, \delta_N f\rangle
+\sum_{l=1}^L\sum_{j<0}\sum_{k\in\N_0}\langle \delta_N f, \tilde\psi^l_{j,k}\rangle \langle\tilde\varphi^l_{j,k}, \delta_N f\rangle.
\end{eqnarray*}
The second term goes to $0$ as $N$ goes to $\infty$, by~(\ref{e.dual}). Hence,
\begin{equation}\label{e.in}
\sum_{l=1}^L\sum_{j\geq 0}\sum_{k\in\N_0}\langle \delta_N f, \tilde\psi^l_{j,k}\rangle \langle\tilde\varphi^l_{j,k}, \delta_N f\rangle\rightarrow \|f\|_2^2~{\rm as}~N\rightarrow\infty.
\end{equation}
But,
\begin{eqnarray*}
\sum_{l=1}^L\sum_{j\geq 0}\sum_{k\in\N_0}\langle \delta_N f, \tilde\psi^l_{j,k}\rangle \langle\tilde\varphi^l_{j,k}, \delta_N f\rangle 
& = & \sum_{l=1}^L\sum_{j\geq 0}\sum_{k\in\N_0}\langle \delta_N f, \psi^l_{j,k}\rangle \langle\varphi^l_{j,k}, \delta_N f\rangle \\
& = & \sum_{l=1}^L\sum_{j\geq 0}\sum_{k\in\N_0}\langle f, \psi^l_{j-N,k}\rangle \langle\varphi^l_{j-N,k}, f\rangle \\
& = & \sum_{l=1}^L\sum_{j\geq -N}\sum_{k\in\N_0}\langle f, \psi^l_{j,k}\rangle \langle\varphi^l_{j,k}, f\rangle.
\end{eqnarray*}
Hence, by~(\ref{e.in}), we have 
\[
\sum_{l=1}^L\sum_{j\in\Z}\sum_{k\in\N_0}\langle f, \psi^l_{j,k}\rangle \langle\varphi^l_{j,k}, f\rangle=\|f\|_2^2.
\]
This shows that~(\ref{e.adual}) holds for all $f=g$ with compact support. Since such functions are dense in $L^2(K)$,~(\ref{e.adual}) holds for all $f=g$ in $L^2(K)$. This completes the proof of the theorem.
\end{proof}
%
%%%%%%%%%%%%%%%%%%%%%%%%%%%%%%%%%%%%%%%%%%%%%%%%%%%%%%%%%%%%%%%%%%%%

\section{Co-affine Systems}
Recall that the quasi-affine system $\tilde X(\Psi)$ was obtained from the affine-system $X(\Psi)$ by reversing the dilation and translation operations for negative scales $j<0$ and then by renormalizing. It is a natural question to ask what happens if we reverse these operations for each scale $j\in\Z$. We make the following definition.

\begin{definition} 
Let $\Psi=\{\psi^1,\psi^2,\dots,\psi^L\}\subset L^2(K)$ and let $c=\{c_{l,j}:1\leq l\leq L, j\in\Z\}$ be a sequence of scalars. The weighted \emph{co-affine system} $X^*(\Psi, c)$ generated by $\Psi$ and $c$ is the collection
\[
X^*(\Psi, c)=\{\psi^{*l}_{j, k}=c_{l,j}\tau_{u(k)}\delta_j \psi^l: 1\leq L, j\in\Z, k\in\N_0\}.
\]
\end{definition}

This concept was first defined by Gressman, Labate, Weiss and Wilson in~\cite{GLWW} for the case of the real line when $\Psi$ consists of a single function and the dilation is a real number greater than 1. They proved that in this case the weighted co-affine system can never be a frame for $L^2(\R)$. Johnson~\cite{J} extended this result to $L^2(\R^n)$ for finitely generated co-affine systems associated with expansive dilation matrices. In this section we will extend this result to the case of a local field of positive characteristic.

Let $X^*(\Psi, c)$ be a weighted co-affine system generated by $\Psi$ and $c$. For $f\in L^2(K)$, define
\[
w_f(x)=\sum_{l=1}^L\sum_{j\in\Z}\sum_{k\in\N_0}|\langle\tau_x f, \psi^{*l}_{j,k}\rangle|^2.
\] 
By Proposition~\ref{p.un}(b) and (c), it follows that $w_f(x+u(n))=w_f(x)$ for all $n\in\N_0$, that is, $w_f$ is integral-periodic. We first prove a lemma. 

\begin{lemma}\label{l.wf}
If $X^*(\Psi, c)$ is a Bessel system for $L^2(K)$, then for each $f\in L^2(K)$, we have
\[
\int_{\D}w_f(x)dx=\int_K\sum_{l=1}^L\sum_{j\in\Z}|c_{l,j}|^2q^{-j}|\hat\psi^l(\p^j\xi)|^2|\hat f(\xi)|^2d\xi.
\]
\end{lemma}

\begin{proof}
The result follows from the Plancherel theorem and the fact that 
\[
(\tau_y\delta_jg)^{\wedge}(\xi)=q^{-j/2}\hat g(\p^j\xi)\overline{\chi_\xi(y)}\quad{\rm for}~y\in K, j\in\Z.
\]
We have,
\begin{eqnarray*}
\int_{\D}w_f(x)dx
& = & \int_{\D}\sum_{l=1}^L\sum_{j\in\Z}\sum_{k\in\N_0}|\langle\tau_x f, c_{l,j}\tau_{u(k)}\delta_j\psi^l\rangle|^2dx \\
& = & \sum_{l=1}^L\sum_{j\in\Z}|c_{l,j}|^2\int_{\D}\sum_{k\in\N_0}|\langle f, \tau_{-x+u(k)}\delta_j\psi^l\rangle|^2dx \\
& = & \sum_{l=1}^L\sum_{j\in\Z}|c_{l,j}|^2\int_K|\langle f, \tau_{-x}\delta_j\psi^l\rangle|^2dx \\
& = & \sum_{l=1}^L\sum_{j\in\Z}|c_{l,j}|^2\int_K\Bigl|\int_K\hat f(\xi)q^{-j/2}\overline{\hat\psi^l(\p^j\xi)}\chi_{x}(\xi)d\xi\Bigr|^2dx \\
& = & \sum_{l=1}^L\sum_{j\in\Z}|c_{l,j}|^2q^{-j}\int_K\Bigl|\Bigl(\hat f~\overline{\hat\psi^l(\p^j\cdot)}\Bigr)^{\vee}(x)\Bigr|^2dx \\
& = & \sum_{l=1}^L\sum_{j\in\Z}|c_{l,j}|^2q^{-j}\int_K|\hat f(\xi)|^2|\hat\psi^l(\p^j\xi)|^2 d\xi.
\end{eqnarray*}
\end{proof}

We now show that there do not exist any co-affine frame in $L^2(K)$.
\begin{theorem}
Let $\Psi=\{\psi^1,\psi^2,\dots,\psi^L\}\subset L^2(K)$ and $c=\{c_{l,j}:1\leq l\leq L, j\in\Z\}$ be a sequence of scalars. Then $X^*(\Psi, c)$ cannot be a frame for $L^2(K)$.
\end{theorem}

\begin{proof}
Suppose that $X^*(\Psi, c)$ is a frame with bounds $A^*$ and $B^*$. That is, 
\[
A^*\|f\|_2^2\leq \sum_{l=1}^L\sum_{j\in\Z}\sum_{k\in\N_0}|\langle f, \psi^{*l}_{j,k}\rangle|^2\leq B^*\|f\|_2^2 \quad~{\rm for~all}~f\in L^2(K).
\]
Taking $f=\psi^{*l_0}_{j_0,k_0}$ for a fixed $j_0\in\Z$, $k_0\in\N_0$ and $1\leq l_0\leq L$, we have
\[
\|\psi^{*l_0}_{j_0,k_0}\|_2^4\leq \sum_{l=1}^L\sum_{j\in\Z}\sum_{k\in\N_0}|\langle \psi^{*l_0}_{j_0,k_0}, \psi^{*l}_{j,k}\rangle|^2\leq B^*\|\psi^{*l_0}_{j_0,k_0}\|_2^2.
\]
This implies $\|\psi^{*l_0}_{j_0,k_0}\|_2^2\leq B^*$. Since $\|\psi^{*l_0}_{j_0,k_0}\|_2=|c_{l_0,j_0}|\|\psi^{l_0}\|_2$, it follows that
\begin{equation}\label{e.clj}
|c_{l_0,j_0}|^2\leq B^*\|\psi^{l_0}\|_2^{-2}\quad~{\rm for~all}~j_0\in\Z~{\rm and}~1\leq l_0\leq L.
\end{equation}
From the definition of $w_f$, we have
\[
A^*\|f\|_2^2\leq w_f(x)\leq B^*\|f\|_2^2\quad~{\rm for~all}~f\in L^2(K).
\]
Integrating over $\D$ and applying Lemma~\ref{l.wf}, we get
\[
A^*\|f\|_2^2\leq\int_K\sum_{l=1}^L\sum_{j\in\Z}|c_{l,j}|^2q^{-j}|\hat\psi^l(\p^j\xi)|^2|\hat f(\xi)|^2 d\xi\leq B^*\|f\|_2^2\quad~{\rm for~all}~f\in L^2(K).
\]
From this we conclude that 
\[
A^*\leq\sum_{l=1}^L\sum_{j\in\Z}|c_{l,j}|^2q^{-j}|\hat\psi^l(\p^j\xi)|^2\leq B^*\quad~{\rm for~a.e.}~\xi\in K.
\]
Now, integrating over $\P^{-1}\setminus\D$ after making the substitution $\xi\rightarrow \p^n\xi, n\in\Z$, we have
\begin{eqnarray*}
A^*|\P^{-1}\setminus\D|
& \leq & \int_{\P^{-1}\setminus\D}\sum_{l=1}^L\sum_{j\in\Z}|c_{l,j}|^2q^{-j}|\hat\psi^l(\p^{j+n}\xi)|^2d\xi \\
&  =   & \int_{\P^{-1}\setminus\D}\sum_{l=1}^L\sum_{j\in\Z}|c_{l,j-n}|^2q^{-j+n}|\hat\psi^l(\p^j\xi)|^2d\xi \\
&  =   & q^{n}\sum_{l=1}^L\sum_{j\in\Z}\int_{\p^j(\P^{-1}\setminus\D)}|c_{l,j-n}|^2|\hat\psi^l(\xi)|^2d\xi.
\end{eqnarray*}
Applying~(\ref{e.clj}), we have
\begin{eqnarray*}
A^*(q-1)
& \leq & q^{n}\sum_{l=1}^L\sum_{j\in\Z}B^*\|\psi^l\|_2^{-2}\int_{\p^j(\P^{-1}\setminus\D)}|\hat\psi^l(\xi)|^2d\xi \\
&   =  & q^{n}\sum_{l=1}^LB^*\|\psi^l\|_2^{-2}\sum_{j\in\Z}\int_{\p^j(\P^{-1}\setminus\D)}|\hat\psi^l(\xi)|^2d\xi \\
&   =  & q^{n}\sum_{l=1}^LB^* = q^{n}LB^*.
\end{eqnarray*}
That is,
\[
A^*(q-1)\leq q^{n}LB^*\quad~{\rm for~each}~n\in\Z.
\]
Letting $n\rightarrow -\infty$, we see that $A^*=0$. Hence, $X^*(\Psi, c)$ cannot be a frame for $L^2(K)$. 
\end{proof}
%
%%%%%%%%%%%%%%%%%%%%%%%%%%%%%%%%%%%%%%%%%%%%%%%%%%%%%%%%%%%%%%%%%%%%

%
%%%%%%%%%%%%%%%%%%%%%%%%%%%%%%%%%%%%%%%%%%%%%%%%%%%%%%%%%%%%%%%%%%%%

\end{document}